\documentclass[fleqn]{mat01}
\usepackage{times,mathtimy,amssymb,latexsym}
\begin{document}

\def\d{\mbox{\rm d}}
\def\e{\mbox{\rm e}}
\def\Sp{\mbox{\rm Sp}}

\newcommand{\bdes}{\begin{enumerate}}
\newcommand{\edes}{\end{enumerate}}
\newcommand{\bit}{\begin{itemize}}
\newcommand{\eit}{\end{itemize}}
\newcommand{\Q}{\mathbb{Q}}
\newcommand{\R}{\mathbb{R}}
\newcommand{\Z}{\mathbb{Z}}
\newcommand{\C}{\mathbb{C}}
\newcommand{\B}{{\mathcal B}}
\newcommand{\E}{{\mathcal E}}
\newcommand{\F}{{\mathcal F}}
\newcommand{\G}{{\mathcal G}}
\newcommand{\I}{{\mathcal I}}
\newcommand{\J}{{\mathcal J}}
\newcommand{\M}{{\mathcal M}}
\newcommand{\Lc}{{\mathcal L}}
\newcommand{\Sm}{{\mathcal S}}
\newcommand{\W}{{\mathcal W}}
\newcommand{\X}{{\mathcal X}}
\newcommand{\id}{\mathrm{Id}}
\newcommand{\mrm}[1]{\mathrm{#1}}
\newcommand{\ov}{\overline}
\newcommand{\ox}{\circ}
\newcommand{\lgra}{\longrightarrow}
\newcommand{\f}{\frac}
\newcommand{\ci}{\subset}
\newcommand{\by}{\times}
\newcommand{\ot}{\otimes}
\newcommand{\op}{\oplus}
\newcommand{\ld}{\ldots}
\newcommand{\cd}{\cdots}
\renewcommand{\a}{\alpha}
\renewcommand{\b}{\beta}
\newcommand{\g}{\gamma}
\newcommand{\gG}{\Gamma}
\newcommand{\ep}{\epsilon}
\newcommand{\Th}{\Theta}
\newcommand{\p}{\pi}
\newcommand{\ta}{\tau}
\newcommand{\s}{\sigma}
\newcommand{\Si}{\Sigma}
\newcommand{\ph}{\phi}
\newcommand{\Ph}{\Phi}
\newcommand{\ps}{\psi}
\newcommand{\Ps}{\Psi}
\newcommand{\om}{\omega}
\newcommand{\Om}{\Omega}
\newcommand{\nb}{\nabla}
\renewcommand{\l}{\lambda}
\renewcommand{\L}{\Lambda}
\newcommand{\what}{\widehat}
\newcommand{\wg}{\wedge}
\newcommand{\iy}{\infty}
\newcommand{\pl}{\partial}
\newcommand{\mb}{\mbox}
\newcommand{\age}{\mathrm{Image \-}}
\newcommand{\D}{\displaystyle}
\newcommand{\vt}[2]{\ensuremath{\D{#1_1, \ld, #1_#2}}}
\newcommand{\pd}[2]{\ensuremath{\D{\f{\pl#1}{\pl#2}}}}
\newcommand{\pe}[2]{\ensuremath{\D{{\pl#1}/{\pl#2}}}}
\mathchardef\gt="313E 
\mathchardef\lt="313C  

\newtheorem{theore}{Theorem}
\renewcommand\thetheore{\arabic{section}.\arabic{theore}}
\newtheorem{theor}[theore]{\bf Theorem}
\newtheorem{lem}[theore]{\it Lemma}
\newtheorem{propo}[theore]{\rm PROPOSITION}
\newtheorem{coro}[theore]{\rm COROLLARY}
\newtheorem{definit}[theore]{\rm DEFINITION}
\newtheorem{probl}[theore]{\it Problem}
\newtheorem{exampl}[theore]{\it Example}
\newtheorem{pot}[theore]{\it Proof of Theorem}

%
%


\setcounter{page}{459} \firstpage{459}

\title{Gromov--Witten invariants and quantum cohomology}

\markboth{Amiya Mukherjee}{Gromov--Witten invariants and quantum
cohomology}

\author{AMIYA MUKHERJEE}

\address{Stat-Math Division, Indian Statistical Institute, 203
B.T. Road, Kolkata~700~108, India\\
\noindent E-mail: amiya@isical.ac.in\\[1.2pc]
\noindent {\it Dedicated to Professor K B Sinha on the occasion of
his 60th \vspace{-1pc}birthday}}

\volume{116}

\mon{November}

\parts{4}

\pubyear{2006}

\Date{}

\begin{abstract}
This article is an elaboration of a talk given at an international
conference on Operator Theory, Quantum Probability, and
Noncommutative Geometry held during December~20--23, 2004, at the
Indian Statistical Institute, Kolkata. The lecture was meant for a
general audience, and also prospective research students, the idea
of the quantum cohomology based on the Gromov--Witten invariants.
Of course there are many important aspects that are not discussed
here.
\end{abstract}


\maketitle

\section{Introduction}

\looseness -1 Quantum cohomology is a new mathematical discipline
influenced by the string theory as a joint venture of physicists
and mathematicians. The notion was first proposed by
Vafa~\cite{V}, and developed by Witten \cite{W} and others
\cite{B}. The theory consists of some new approaches to the
problem of constructing invariants of compact symplectic manifolds
and algebraic varieties. The approaches are related to the ideas
of a $(1+1)$-dimensional topological quantum field theory, which
indicate that the general principle of constructing invariants
should be as follows: The invariants of a manifold $M$ should be
obtained by integrating cohomology classes over certain moduli
space $\M$ associated to $M$.  In our case the manifold is a
symplectic manifold $(M,\om)$, and the moduli space $\M$ is the
space of certain $J$-holomorphic spheres $\s\hbox{:}\ \C P^1\lgra
M$ in a given homology class $A\in H_2(M;\Z)$.  The relevant
cohomology classes on $\M$ are the pullbacks $e^* a$ of the
cohomology classes $a\in H^*(M,\Z)$ under the evaluation maps
$e\hbox{:}\ \M\lgra M$ given by $e(\s)=\s(z)$ for fixed $z\in \C
P^1$. Then the integration of a top dimensional product of such
classes (or equivalently, the evaluation of the top dimensional
form on the fundamental class $[\M]$) gives rise to the
Gromov--Witten invariant
\begin{equation*}
\int_\M e^*a_1 \wg \cd \wg e^* a_p = \langle e^*a_1 \wg \cd \wg
e^* a_p, [\M]\rangle.
\end{equation*}
These invariants are independent of the choices of $\vt{z}{p}$ in
$\C P^1$ used in their definitions, and can be interpreted as
homomorphisms
\begin{equation*}
\Ph_A\hbox{:}\ H_*(M,\Z)\ot \cd \ot H_*(M,\Z)\lgra \C
\end{equation*}
given by
\begin{equation*}
\Ph_A(\vt{\a}{p})=\int_\M e^* a_1 \wg \cd \wg e^*  a_p,
\end{equation*}
where $a_j$ is the $\mrm{Poincar\acute{e}}$ dual of
$\a_j$.\pagebreak

The invariant $\Ph_A$ counts the number of intersection points
(with signs of their orientations) of the image of the $p$-fold
evaluation map $\s\mapsto (\s(z_1), \ld , \s(z_p))\in M^p$ with
the cycles representing the $\a_j$, where the dimension of the
homology class
\begin{equation*}
\a_1\by \cd \by \a_p
\end{equation*}
is chosen so that if all the intersections were transversal, there
would be only a finite number of such points. This is simply the
number of $J$-holomorphic spheres in the given homology class $A$
which meets the cycles representing $\vt{\a}{p}$.

The importance of the $J$-holomorphic spheres and the
Gromov--Witten invariants is that they may be used to define a
quantum deformation of the cup product in the cohomology ring
$H^*(M)$ of a compact symplectic manifold $M$ making it a quantum
cohomology ring $QH^*(M)$.

The description of the Gromov--Witten invariants can be given in
terms of a general Riemann surface $\Si$ (see~\cite{RT}). However,
we have made this expository introduction somewhat simpler by
taking $\Si=S^2$. The results that guided our approach are to be
seen in the work of MacDuff and Salamon~\cite{MS2}.

\section{\pmb{$J$}-Holomorphic curves}

A symplectic manifold $(M,\om)$ is a smooth manifold $M$ of
dimension $2n$ with a symplectic structure $\om$ on it which is a
closed differential $2$-form $\om$ such that the volume form
$\om^n$ is nowhere vanishing on $M$. The basic example is the
Euclidean space $\R^{2n}$ with the constant symplectic form
\begin{align*}
\om_0= \d x_1\wg \d y_1 + \d x_2\wg \d y_2 + \cd + \d x_n\wg \d
y_n,
\end{align*}
where $(\vt{x}{n},\vt{y}{n})$ are coordinates in $\R^{2n}$. The
next basic example is provided by the phase space of a Hamiltonian
system, that is, by the cotangent bundle $T^*N$ of any
$n$-manifold $N$ with a symplectic structure which is locally the
pullback of the structure $\om_0$ on $\R^{2n}$. A~sympletic
manifold cannot be odd dimensional.

A symplectic diffeomorphism $\ph\hbox{:}\ (M_1,\om_1)\lgra
(M_2,\om_2)$ between two symplectic manifolds is a diffeomorphism
$\ph\hbox{:}\ M_1\lgra M_2$ such that $\ph^*\om_2=\om_1$.
Symplectic geometry is quite different from Riemannian geometry,
and also from K\"{a}hlerian geometry. The Darboux theorem says
that locally any two symplectic manifolds of the same dimension
are diffeomorphic.  Therefore locally all symplectic manifolds are
the same, and there is no local invariant in symplectic geometry,
like, for example, the curvature in Riemannian geometry. The only
possible invariants have to be global.

The Darboux theorem makes it difficult to study the global
structure on a symplectic manifold. Although variational
techniques may be employed to tackle some global problems, it is
the theory of $J$-holomorphic curves of Gromov that applies to
many problems of symplectic manifolds. We have no other theory to
investigate these global\break questions.

An almost complex structure on a manifold $M$ is a complex
structure on its tangent bundle $TM$, that is, an endomorphism
$J\hbox{:}\ TM\lgra TM$ such that $J^2=-\id$. Then, $J$ makes $TM$
a complex vector bundle, where the complex vector space structure
on each fibre $T_xM$ is given by $(a+\sqrt{-1}b)\cdot v=a v +bJ
v$. If $M$ is already a complex manifold, which is a manifold with
holomorphic changes of coordinates, then the tangent bundle $TM$
is a complex vector bundle, and its almost complex structure $J$
is just multiplication by $\sqrt{-1}$. The standard almost complex
structure $j$ on $\R^{2n}$ is given by
\begin{align*}
j(\pe{}{x_k})=\pe{}{y_k},\ \quad\ j(\pe{}{y_k})= -\pe{}{x_k},
\end{align*}
where $z_k=x_k+\sqrt{-1} y_k$ are the coordinates in $\C^n$.

An almost complex structure $J$ on a symplectic manifold $(M,\om)$
is called $\om$-tamed if $\om(v, Jv) > 0$ for every nonzero vector
$v\in T_xM$. This means that the restriction of $\om$ to the
complex line in $T_xM$ spanned by $v$ and $Jv$ verifies the same
condition, and so $\om$ restricts to a non-degenerate form on each
such line. An almost complex structure $J$ is called
$\om$-compatible if it is $\om$-tamed, and
\begin{align*}
\om(Jv,Jw)=\om(v,w)\quad \mb{ for all} \ v, w\in TM.
\end{align*}

The space $\J(M,\om)$ of all $\om$-compatible almost complex
structures on $(M,\om)$ is non-empty and contractible, because
associated with the tangent bundle $TM$ we have a bundle
$\J(M,\om)\lgra M$ with contractible fibre $\Sp(2n,\R)/U(n)$.
Since $\J(M,\om)$ is pathwise connected, the complex bundles
$(TM,J)$ are isomorphic for different choices of $J\in \J(M,\om)$.
Therefore the Chern classes $c_i(M)$ of these bundles do not
depend on $J$. The assertions also apply to $\om$-tamed almost
complex structures (in this case the associated bundle has fibre
$GL(2n,\R)/U(n)$).


A smooth map $\ph\hbox{:}\ (M,J)\lgra (M',J')$ between almost
complex manifolds is called $(J,J')$-holomorphic if
$\d\ph_x\hbox{:}\ T_xM\lgra T_{\ph(x)}M'$ is complex linear, that
is, $\d\ph_x\ox J_x= J'_{\ph(x)}\ox \d\ph_x$ for all $x\in M$.
These conditions are exactly the Cauchy--Riemann equations in the
case when $(M,J)$ and $(M',J')$ are subsets of $\C^n$. An almost
complex structure $J$ on $M$ is called integrable if it arises
from a complex structure on $M$; in other words, if $M$ admits an
atlas whose coordinate charts are $(J,j)$-holomorphic maps, where
$j$ denotes the standard complex structure on $\C^n$.  If $\dim
M=2$, a fundamental theorem says that any almost complex structure
$J$ on $M$ is integrable. However, the theorem is not true in
higher dimensions. The non-integrability of $J$ is measured by the
Nijenhuis tensor $N_J$ (see~\cite{MS1}).

A $J$-holomorphic curve in $(M,J)$ is a $(J_0,J)$-holomorphic map
$\s\hbox{:}\ \Si\lgra M$, where $(\Si,J_0)$ is a Riemann surface
(complex manifold of dimension $1$) with complex structure $J_0$.
Very often we take $(\Si,J_0)$ as the Riemann sphere $S^2$, and in
this case a $J$-holomorphic curve is referred to as a
$J$-holomorphic sphere.

If $\s$ is an embedding and $C$ is the image of $\s$, then $\s$ is
called a $J$-holomorphic parametrization of $C$. In this case $C$
is $2$-submanifold of $M$ with $J$-invariant tangent bundle $TC$
so that each tangent space is a complex line in $TM$. Conversely,
any $2$-submanifold $C$ of $M$ with a $J$-invariant tangent bundle
has a $J$-holomorphic parametrization $\s$, because the
restriction of $J$ to $C$ is integrable.

For an $\om$-tamed almost complex structure $J$ on a symplectic
manifold $(M,\om)$, the image of $J$-holomorphic parametrization
is a symplectic $2$-submanifold of $M$ with $J$-invariant tangent
spaces. Conversely, given an oriented $2$-submanifold $C$ of $M$,
one can construct an $\om$-tamed $J$ such that $TC$ is
$J$-invariant (first define $J$ on $TC$ and then extend it to TM).
One may contrast this situation with that in complex geometry
where one often defines a curve as the set of common zeros of a
number of holomorphic polynomials. Such an approach makes no sense
in the case when the almost complex structure is non-integrable,
since there may not exist holomorphic functions $(M,J)\lgra\C$
when $J$ is non-integrable.

\section{Moduli spaces}

Let $(M,J)$ be an almost complex manifold without boundary, and
$\Si$ be a Riemann surface of genus $g$ with complex structure
$J_0$. Then a moduli space $\M(A,J)$ is the space of all simple
$J$-holomorphic curves $\s\hbox{:}\ (\Si,J_0)\lgra (M,J)$ which
represent a given homology class $A\in H_2(M;\Z)$ (i.e.
$\s_*[\Si]=A$), with the $C^r$-topology, $r\ge 0$. Our first
problem is to provide a finite dimensional smooth structure on
this space.

Note that a curve $\s$ is simple, if it is not multiply-covered,
that is, it is not a composition of a holomorphic branched
covering $(\Si,J_0)\lgra (\Si',J'_0)$ of degree $\gt\ 1$ and a
$J$-holomorphic map $\Si'\lgra M$. We avoid multiply-covered
curves because they may be singular points in the moduli space
$\M(A,J)$. Every simple curve $\s$ has an injective point $z\in
\Si$, which is a regular point of $\s$ (i.e. $\d\s_z\ne 0$) such
that $\s^{-1}\s(z)=\{z\}$. Moreover, the set of injective points
is open and dense in $\Si$~\cite{MS2}.

The space $\Sm=C^\iy(\Si,M,A)$ of smooth maps $\s\hbox{:}\
\Si\lgra M$ that are somewhere injective, and represent $A\in
H_2(M;\Z)$ may be looked upon as an infinite dimensional manifold
whose tangent space at $\s\in \Sm$ is given by
\begin{equation*}
T_\s\Sm=C^\iy(\s^*TM),
\end{equation*}
which is the vector space of all smooth vector fields of $M$ along $\s$.

We can view $\s^*TM$ as a complex vector bundle. Therefore we have
a splitting of the space of $1$-forms
\begin{equation*}
\Om^1(\s^*TM)=\Om^{1,0}(\s^*TM)\op \Om^{0,1}(\s^*TM),
\end{equation*}
where $\Om^{1,0}$ and $\Om^{0,1}$ are respectively vector spaces
of $J$-linear and $J$-anti-linear $1$-forms with values in
$\s^*TM$. Since $\d\s\in\Om^1(\s^*TM)$, we can decompose
$\d\s=\pl_J(\s) + \ov{\pl}_J(\s)$,~where
\begin{align*}
\pl_J(\s) &= \f{1}{2}(\d\s - J\ox \d\s\ox J_0),\\[.4pc]
\ov{\pl}_J(\s) &= \f{1}{2}(\d\s + J\ox \d\s\ox J_0)
\end{align*}
are respectively $J$-linear and $J$-anti-linear parts of $\d\s$.

If $\E\lgra \Sm$ is the infinite dimensional vector bundle whose
fibre $\E_\s$ over $\s\in\Sm$ is the space $\Om^{0,1}(\s^*TM)$,
then $\ov{\pl}_J$ is a section of the bundle $\E\lgra\Sm$.
Moreover, the $J$-holomorphic curves are the zeros of the section
$\ov{\pl}_J$, that is, if $Z$ denotes the zero section of the
bundle, then
\begin{equation*}
\M(A,J)=(\ov{\pl}_J)^{-1}(Z).
\end{equation*}
This will be a manifold if $\ov{\pl}_J\hbox{:}\ \Sm\lgra\E$ is
transversal to $Z$, that is, the image of
\begin{equation*}
\d\ov{\pl}_J(\s)\hbox{:}\ T_\s \Sm \lgra T_{(\s,0)}\E
\end{equation*}
is complementary to the tangent space of the zero-section $Z$ for
every $\s\in \M(A,J)$; in other words, the linear operator $D_\s=
\p_\s\ox \d\ov{\pl}_J(\s)$, where
\begin{equation*}
\p_\s\hbox{:}\ T_{(\s,0)}\E = T_\s \Sm\op \E_\s\lgra\E_\s
\end{equation*}
is the projection, is surjective for every
$\s\in\M(A,J)$.\pagebreak

Explicit expression of the operator
\begin{equation*}
D_\s\hbox{:}\ C^\iy(\s^* TM)\lgra \Om^{0,1}(\s^* TM),
\end{equation*}
can be obtained by differentiating the local expressions of
$\ov{\pl}_J(\s)$ in the direction of a vector field along $\s$.
These expressions show that the first order terms make up the
usual Cauchy--Riemann operator for maps $\C\lgra \C^n=\R^{2n}$.
Therefore $D_\s$ is a first-order elliptic differential operator,
and hence it is Fredholm. Recall that a bounded operator
$F\hbox{:}\ X\lgra Y$ between Banach spaces $X$ and $Y$ is a
Fredholm operator if $F$ has finite dimensional kernel and
cokernel, and $F(X)$ is closed. The index of $F$ is defined by
\begin{equation*}
\mrm{index}\,  F = \dim\,\ker F - \dim\  \mrm{coker}\, F.
\end{equation*}
These operators form an open subset $\F(X,Y)$ of the space of
bounded operators $\B(X,Y)$ with the norm topology. The Fredholm
index is constant on each connected component of $\F(X,Y)$, and
therefore $\mrm{index}\, F$ is not altered if $F$ varies
continuously.

Although the domain and range of the Fredholm operator $D_\s$ are
complex vector spaces, $D_\s$ is not complex linear, because $J$
is not integrable. It will appear from the computations for $D_\s$
mentioned above that the complex anti-linear part of $D_\s$ has
order $0$. Then, by multiplying the anti-linear part by a constant
which tends to $0$, we can find a homotopy of $D_\s$ through
Fredholm operators. The final Fredholm operator of the homotopy
commutes with $J$, and is a Cauchy--Riemannian operator. It
determines a holomorphic structure on the complex vector bundle
$\s^*TM$. Therefore we have by the Riemann--Roch theorem
(\cite{GH}, p.~243)
\begin{equation*}
\mrm{index}\ D_\s = n(2-2g) + 2 c_1(\s^* TM)[\Si] = n(2-2g) + 2
c_1(A),
\end{equation*}
where $c_1$ is the first Chern class of the complex bundle
$(TM,J)$, and $c_1(\s^*TM)[\Si] = (\s^*c_1) [\Si] =
c_1(\s_*[\Si])=c_1(A)$.

If the operator $D_\s$ is surjective for every $\s\in \M(A,J)$,
then it follows from the infinite dimensional implicit function
theorem that $\M(A,J)$ is a finite dimensional manifold whose
tangent space at $\s$ is $\ker D_\s$.

We suppose that the space of $\om$-compatible almost complex
structures $\J=\J(M,\om)$ has been endowed with the
$C^\iy$-topology.  Let $\J_r$ be the subspace of $\J$ consisting
of those structures $J$ for which $D_\s$ is surjective for all
$\s\in \M(A,J)$.

\begin{theor}[\!] $\left.\right.$
\begin{enumerate}
\leftskip .12pc
\renewcommand\labelenumi{{\rm (\alph{enumi})}}
\item If $J\in \J_r${\rm ,} then $\M(A,J)$ is a smooth manifold with a
natural orientation such that
\begin{equation*}
\dim \M(A,J)= n(2-2g)+2 c_1(A).
\end{equation*}

\item The subset $\J_r$ is residual in $\J$.
\end{enumerate}
\end{theor}

Recall that a subset of a topological space $X$ is residual if it
is the intersection of a countable family of open dense subsets of
$X$. A~point of $X$ is called generic if it belongs to some
residual subset of $X$.

\begin{proof}
Part~(a) follows from the above discussion, except for the
orientation. The orientation follows from the fact that a Fredholm
operator $D$ between complex Banach spaces induces a canonical
orientation on its determinant line
\begin{equation*}
\det\, D=\Om^p(\mrm{ker}\ D)\ot \Om^q(\mrm{ker}\ D^*),
\end{equation*}
where $p=\dim \mrm{ker}\ D$ and $q=\dim \mrm{coker}\ D$, provided
$D$ is complex linear. As described above, we may suppose by using
a homotopy of order $0$ that our Fredholm operator $D_\s$ is
complex linear with $\mrm{coker}\ D_\s=0$. Therefore its
determinant line, and hence $\mrm{ker}\ D_\s = T_\s \M(A,J)$, has
a canonical orientation.  These arguments are due to
Ruan~\cite{R}, also note that earlier Donaldson~\cite{D} used
similar arguments for the orientation of \hbox{Yang--Mills} moduli
spaces.

Part~(b) uses an infinite dimensional version of Sard--Smale
theorem which is due to Smale~\cite{S}. A~non-linear smooth map
$f\hbox{:}\ X\lgra Y$ between Banach spaces is a Fredholm map of
index $k$, if the derivative $\d f_x\hbox{:}\ X\lgra Y$ is a
linear Fredholm operator of index $k$ for each $x\in X$.  A point
$y\in Y$ is a regular value of $f$ if $\d f_x$ is surjective for
each $x\in f^{-1}(y)$, otherwise $y$ is called a critical value of
$f$. Then the Sard--Smale theorem says that if $f\hbox{:}\ X\lgra
Y$ is a $C^k$ Fredholm map between seperable Banach spaces and
$k\gt \max(0, \mrm{index} \, f)$, then the set of regular values
of $f$ is residual in $Y$. The theorem remains true if $X$ and $Y$
are Banach manifolds, instead of Banach spaces. It follows from
the implicit function theorem for Banach spaces that if $y\in Y$
is a regular value then $f^{-1}(y)$ is a smooth submanifold of
$X$. Moreover, if $f^{-1}(y)$ is finite dimensional, then its
dimension is equal to the Fredholm index of $f$.

For the completion of the proof of the theorem, we need to refine
the space $\Sm$ using the Sobolev $W^{k,p}$-norm which is given by
the sum of the $L^p$-norms of all derivatives of $\s\in \Sm$ up to
order $k$,
\begin{equation*}
\|\s\|_{k,p}=\sum_{|r|\le k}\  L^p(\pl^r \s),
\end{equation*}
where $r$ is a multi-index and $|r|$ is its order. It can be shown
that the Sobolev space $W^{k,p}(\Si,M)$, which is the space
consisting of all maps $\Si\lgra M$ whose $k$-th order derivatives
are of class $L^p$ (and which represent the class $A\in
H_2(M;\Z)$), is the completion of the space $\Sm$ with respect to
the Sobolev $W^{k,p}$-norm (see Appendix~B in~\cite{MS2}). It
appears that we must assume the condition $kp\gt 2$ in order for
the space $W^{k,p}(\Si,M)$ to be well-defined. Under this
condition, the Sobolev embedding theorem says that there is a
continuous embedding of $W^{k,p}(\Si,M)$ into the space of
continuous maps $C^0(\Si,M)$, and the multiplication theorem says
that the product of two maps of class $W^{k,p}$ is again a map of
the same class.

At the same time we restrict the space of almost complex
structures $\J(M,\om)$, introduced earlier. Let $\J^\ell$,
$\ell\ge 1$, be the space of all almost complex structures of
class $C^\ell$ which are compatible with $\om$, with the $C^\ell$
topology. We shall choose $\ell$ later according to our
requirement.

Then $\J^\ell$ is a smooth separable Banach manifold. Let ${\rm
End}(TM,J,\om)\lgra M$ be the bundle whose fibre over $p\in M$ is
the space of linear endomorphisms $X\hbox{:}\ T_pM\lgra T_pM$ such
that
\begin{equation*}
XJ + JX=0, \quad \om(Xv,w) +\om(v,Xw)=0,\quad \mbox{ for } \ \ \
v, w\in T_pM.
\end{equation*}
Then the tangent space $T_J \J^\ell$ at $J$ is the space of
sections of this bundle.

It can be proved by elliptic bootstrapping methods
(see~\cite{MS2}, Appendix~B for details) that if $J\in \J^\ell$
with $\ell\ge 1$, then a $J$-holomorphic curve $\s\hbox{:}\
\Si\lgra M$ of class $W^{\ell,p}$ with $p\gt 2$ is also of class
$W^{\ell+1,p}$. In particular, if $J$ is smooth and $\s$ is of
class $C^\ell$, then $\s$ is also smooth. Thus if $k\le \ell+1$
and $J\in\J^\ell$, then the moduli space of $J$-holomorphic curves
of class $W^{k,p}$ does not depend on $k$.

In the context of the Sobolev space of $W^{k,p}$-maps $\s\hbox{:}\
\Si\lgra M$ for some fixed $p\gt 2$, we have the Banach space
bundle $\E^p\lgra W^{k,p}(\Si,M)$ whose fibre over $\s\in
W^{k,p}(\Si,M)$ is the space
\begin{equation*}
\E^p_\s = L^p(\L^{0,1}T^*\Si\ot_J \s^* TM)
\end{equation*}
of complex anti-linear $1$-forms on $\Si$ of class $L^p$ taking
values in $\s^*TM$. The non-linear Cauchy--Riemann equations
determine a section $\ov{\pl}_J$ of this bundle, and the
derivative of $\ov{\pl}_J$ at $\s$ gives rise to the operator
\begin{equation*}
D_\s\hbox{:}\  W^{k,p}(\s^*TM)\lgra W^{k-1,p}(\L^{0,1}T^*\Si\ot_J
\s^* TM).
\end{equation*}
The explicit formula for $D_\s$ is given by
\begin{equation*}
D_\s\xi =\f{1}{2}(\nb\xi + J(\s)\nb\xi\ox J_0) + \f{1}{8}
N_J(\pl_J(\s), \xi),
\end{equation*}
where $\nb$ is the Hermitian connection on $M$, and $N_J$ is the
Nijenhuis tensor (see~\cite{M1}). The first part has order 1 and
commutes with $J$, while the second has order 0 and anti-commutes
with $J$.

The ellipticity of $D_\s$ can be established from the estimate
\begin{equation*}
\| \xi \|_{W^{1,p}}\le c_0 (\|D_\s\xi\|_{L^p}+\|\xi\|_{L^p}),
\end{equation*}
which follows from the $L^p$-estimate for Laplace operator (the
Calderson--Zygmund inequality) (see Appendix~B in~\cite{MS2}).
Therefore $D_\s$ is a Fredholm operator of positive index, by a
previous argument in a similar situation.

The following space is also a smooth Banach manifold
\begin{equation*}
\M^\ell(A,\J^\ell)=\{(\s,J) \in W^{k,p}(\Si,M)\by
\J^\ell|\ov{\pl}_J(\s)=0\}.
\end{equation*}
The tangent space $T_{(\s,J)}\ \M^\ell(A,\J^\ell)$ is the space of
all pairs $(X,Y)$ such that
\begin{equation*}
D_\s X + \f{1}{2}\ Y(\s)\ox \d\s\ox J_0=0.
\end{equation*}

Let $\p\hbox{:}\ \M^\ell(A,\J^\ell)\lgra \J^\ell$ be the
projection. Then $\p^{-1}(J)= \M^\ell(A,J)$, and the derivative of
$\p$ at $(\s,J)$,
\begin{equation*}
\d\p(\s,J)\hbox{:}\ T_{(\s,J)}\M^\ell(A,\J^\ell)\lgra T_J\J^\ell
\end{equation*}
is just the projection $(X,Y)\mapsto Y$. It follows that
$\d\p(\s,J)$ is a Fredholm operator having the same index as
$D_\s$. Moreover, a regular value $J$ of $\p$ is an almost complex
structure such that $D_\s$ is surjective for all $J$-holomorphic
spheres $\s\in \p^{-1}(J)$.

We denote the set of regular values of $\p$ by $\J^\ell_r$. By the
Sard--Smale theorem (stated earlier), the set $\J^\ell_r$ is
residual in $\J^\ell$ with respect to the $C^\ell$ topology
whenever $\ell-2 \ge \mrm{index}\ D_\s = \mrm{index}\  \p$,
because $\p$ is of class $C^{\ell -1}$.

Let $\l$ be a positive number, and $\J^\ell_{r,\l}$ be the set of
almost complex structures $J\in \J=\J(M,\om)$ such that $D_\s$ is
surjective for every $J$-holomorphic sphere $\s$ with
$\|\s\|_{L^\iy}\le \l$. Clearly, the intersection of the sets
$\J^\ell_{r,\l}$ over all $\l\gt 0$ is the set $\J_r$ of part~(b)
of the theorem.

The set $\J_r$ is residual, because each $\J^\ell_{r,\l}$ is open
and dense in $\J$ with respect to the $C^\iy$ topology. We omit
the details which may be found in \S 3.4 of~\cite{MS2}. \hfill
$\Box$
\end{proof}

The above theorem can be extended further in order to understand
how the manifolds $\M(A,J)$ depend on $J\in \J_r$. Two almost
complex structures $J_0$ and $J_1$ in $\J_r$ are called smoothly
homotopic if there is a smooth path $[0,1]\lgra \J$, $t\mapsto
J_t$, from $J_0$ to $J_1$.

\begin{theor}[\!]
Let $\J$ be path-connected, and $J_0${\rm ,} $J_1\in \J_r$. Let
$\J(J_0, J_1)$ be the space of all smooth homotopies from $J_0$ to
$J_1$. Then there is a dense set
\begin{equation*}
\J_r(J_0, J_1)\ci \J(J_0,J_1)
\end{equation*}
such that for every $\{J_t\}\in \J_r(J_0,J_1)${\rm ,} the space
\begin{equation*}
\M(A,\{J_t\}_{t\in [0,1]})=\{(t,\s)|\s\in \M(A,J_t)\}
\end{equation*}
is a smooth manifold of dimension $n(2-2g)+ 2c_1(A) +1$ with a
natural orientation and with a smooth boundary which is given by
\begin{equation*}
\pl\M(A, \{J_t\}_{t\in [0,1]})=\M(A,J_1) - \M(A,J_0),
\end{equation*}
where the negative sign indicates the reversed orientation.
\end{theor}

Thus moduli spaces $\M(A,J_0)$ and $\M(A,J_1)$ are oriented
cobordant. We may call the elements of the set $\J_r(J_0, J_!)$
regular homotopies.

\section{Compactness}

The manifold $\M(A,J)$ will not serve any purpose unless some kind
of compactness is established for it.

For simplification we suppose that $\M(A,J)$ is the moduli space
of $J$-holomorphic spheres. The Rellich's theorem says that the
inclusion map
\begin{equation*}
W^{k+1,p}(S^2,M)\lgra W^{k,p}(S^2,M)
\end{equation*}
is compact for all $k$ and $p$ (this means that a sequence
$\{\s_n\}$ which is bounded in the domain $W^{k+1,p}(S^2,M)$
possesses a subsequence which is convergent in the range
$W^{k,p}(S^2,M)$). Moreover, if $k-2/p\gt m+\a$ where $0\lt \a\lt
1$, then $W^{k,p}(S^2,M)$ embeds compactly into the
H$\ddot{\rm{o}}$lder space $C^{m+\a}(S^2,M)$. Using this one gets
the main elliptic regularity theorem, which contains a result of
compactness.

\setcounter{theore}{0}
\begin{theor}[\!]
If $k\ge 1${\rm ,} $p\gt 2${\rm ,} and $\s\in W^{k,p}(S^2,M)$ with
$\ov{\pl}_J\s=0${\rm ,} then $\s\in C^\iy(S^2,M)$. Moreover{\rm ,}
for every integer $m\gt 0${\rm ,} every subset of
${\ov{\pl}_J}^{-1}(0)$ which is bounded in $W^{k,p}(S^2,M)$ has
compact closure in $C^m(S^2,M)$.
\end{theor}

The details are in~\cite{M2}.

An $\om$-tamed almost complex structure $J$ determines a
Riemannian metric on $M$,
\begin{equation*}
\langle v,w\rangle_{J} =\f{1}{2}[\om(v,Jw) +\om(w,Jv)].
\end{equation*}
The energy of a $J$-holomorphic sphere $\s\hbox{:}\ S^2\lgra M$
with respect to this metric is
\begin{equation*}
E(\s)=\int_{S^2}\ |\d\s|^2_J.
\end{equation*}
The group $G=PSL(2,\C)$ acts on $\C\cup\{\iy\}=\C_\iy$ by
$\mrm{M\ddot{o}bius}$ transformations $\ph_L\hbox{:}\
\C_\iy\lgra\C_\iy$,
\begin{equation*}
\ph_L(z)=\f{az+b}{cz+d},\quad\ L=\begin{pmatrix} a & b \\[.1pc] c & d \
\end{pmatrix}  \in SL(2,\C).
\end{equation*}
We may identify $S^2$ with $\C P^1\simeq \C_\iy$ by a
stereographic projection $\p$, and different choices of $\p$
correspond to the action of $SO(3)\simeq SU(2)/\{\pm\id\}\ci
PSL(2,\C)=G$ on $\C P^1=\C_\iy$. Then a $J$-holomorphic sphere
$S^2\lgra M$ gets identified with a smooth $J$-holomorphic curve
$\s\hbox{:}\ \C\lgra M$ such that the map $\C-\{0\}\lgra M$ given
by $z\mapsto \s(1/z)$ extends to a smooth map $\C\lgra M$. The
space of such maps remain invariant under composition with
$\mrm{M\ddot{o}bius}$ transformations $\ph_L\hbox{:}\
\C_\iy\lgra\C_\iy$. We say that a sequence of such $J$-holomorphic
curves $\s_n\hbox{:}\ \C\lgra M$ converges on $\C_\iy$ if both the
sequences $\{\s_n(z)\}$ and $\{\s_n(1/z)\}$ converge uniformly
with all derivatives on compact subsets of $\C$.

It can be shown that
\begin{equation*}
E(\s)=\int_\C\ \s^* \om=\om(A)
\end{equation*}
for all $J$-holomorphic curves $\s\hbox{:}\ \C\lgra M$ ($J$ is
$\om$-tamed), where $\om$ is considered as an integral valued
form. Thus the $L^2$-norm of the derivative of $\s$ satisfies a
uniform bound which depends only on the homology class $A$
represented by $\s$. This does not imply compactness of the moduli
space by the Sobolev estimate, because here $p=2$ (a uniform bound
on the $L^p$-norms of $\d\s$ with $p\gt 2$ would guarantee the
compactness).

It may be noted that the space $\M(A,J)$ can never be both compact
and non-empty, unless $A=0$ in which case all $\s$ are constant
maps. Because, the group $G=PSL(2,\C)$ of holomorphic maps
$S^2\lgra S^2$ is non-compact and it acts on $S^2$ by
reparametrization $\s\mapsto \s\ox\ph$, $\ph\in G$, and so any
$\s\in \M(A,J)$ has a non-compact orbit. However, it is possible
to compactify the quotient $\M(A,J)/G$ sometimes, if $A$ satisfies
a certain condition.

One can show that if $\s_n$ is a sequence in $\M(A,J)$ without any
limit point in $\M(A,J)$, then there is a point $z\in S^2$ such
that the derivatives $\d\s_n(z)$ are unbounded. This implies after
passing to a subsequence that there is a decreasing sequence of
neighbourhoods $U_n$ of $z$ in $S^2$ such that the images
$\s_n(U_n)$ converge to a $J$-holomorphic sphere. If $B$ is the
homology class of this sphere, then either $\om(B)=\om(A)$, or
else $0\lt \om(B)\lt \om(A)$. In the first case, the maps $\s_n$
can be reparametrized so that they converge in $\M(A,J)$. The
second case is referred to as the phenomenon of `bubbling off'.
Here one must proceed with more care. The phenomenon was
discovered by Sacks and Uhlenbeck~\cite{SU} in the context of
minimal surfaces.

The following theorem gives a criterion for the moduli space
$\M(A,J)/G$ to be compact. This is the simplest version of
Gromov's compactness theorem.

A homology class $B\in H_2(M,\Z)$ is called spherical if it lies
in the image of the Hurewicz homomorphism $\p_2(M)\lgra
H_2(M,\Z)$. It is customary to write $B\in\p_2(M)$ if $B$ is a
spherical homology class.

\begin{theor}[\!]
If there is no spherical homology class $B\in H_2(M;\Z)$ such that
$0\lt \om(B)\lt \om(A)${\rm ,} then the moduli space $\M(A,J)/G$
is compact.
\end{theor}

The proof consists of showing that if $\s_n\hbox{:}\ \C\cup
\{\iy\}\lgra M$ is a sequence of $J$-holomorphic $A$-spheres, then
there is a sequence of matrices $L_n\in SL(2,\C)$ such that the
sequence $\s_n\ox\ph_{L_n}$ has a convergent subsequence.
Therefore if $\om(A)$ is already the smallest positive value taken
by $\om$, then the moduli space is compact.

If the criterion of the theorem is not satisfied, it is still
possible sometimes to compactify $\M(A,J)/G$ by adding suitable
pieces. This we shall explain in the next section in a more
general context.

\section{Evaluation maps}

The Gromov--Witten invariants are constructed from the evaluation
map
\begin{equation*}
\M(A,J)\by S^2\lgra M
\end{equation*}
given by $(\s,z)\mapsto \s(z)$. The group $G=PSL(2,\C)$ acts on
the space $\M(A,J)\by S^2$ by $\ph\cdot
(\s,z)=(\s\ox\ph^{-1},\ph(z))$. Therefore we get a map by passing
to the quotient
\begin{equation*}
e=e_J\hbox{:}\ \W(A,J)=\M(A,J)\by_G S^2\lgra M.
\end{equation*}

For example, suppose that $M=\C P^1\by V$ with a product
symplectic form, and $A=[\C P^1\by \{\rm{point}\}]$. If
$\p_2(V)=0$, then $A$ generates a spherical $2$-class in $M$, and
so $\om(A)$ is necessarily the smallest value assumed by $\om$ on
the spherical classes. Therefore by Theorems~3.1 and 4.2, the
space $\W(A,J)$ is a compact manifold for generic $J$. Since
$c_1(A)=2$, $\dim \W(A,J)=2n$ which is the dimension of $M$. It
can be shown that different choices of $J$ give rise to cobordant
maps $e_J$. Since the cobordant maps have the same degree, $\deg
e_J$ is independent of all choices. In the case when $J=J_0\by J'$
is a product, where $J_0$ is the standard complex structure on $\C
P^1$, it can be seen that the elements of $\M(A,J)$ have the form
$\s(z)=(\ph(z),v_0)$, where $v_0\in V$ and $\ph\in G$. It follows
that the map $e_J$ has degree $1$ for this choice of $J$ and hence
for every $J$.

In general, we have a $p$-fold evaluation map
\begin{equation*}
e_p\hbox{:}\ \W(A,J,p)=\M(A,J)\by_G (\C P^1)^p\lgra M^p
\end{equation*}
defined by
\begin{equation*}
e_p(\s, \vt{z}{p})=(\s(z_1), \ld, \s(z_p)).
\end{equation*}
Here, for a space $X$, $X^p$ denotes the $p$-fold
product $X\by \cd\by X$.

For a generic almost complex structure $J$, the space $\W(A,J,p)$
is a manifold with
\begin{equation*}
\dim\ \W(A,J,p)= 2n + 2c_1(A) + 2p -6.
\end{equation*}
This manifold is not compact in general. However, in many cases
the image
\begin{equation*}
\X(A,J,p)=e_p(\W(A,J,p))\ci M^p
\end{equation*}
can be compactified by adding suitable pieces of dimensions at
most equal to $\dim\ \W(A,J,p) -2$. These pieces are called
cusp-curves (the terminology is due to Gromov~\cite{G}), and they
are connected unions of certain $J$-holomorphic spheres. By the
Gromov compactness theorem (which is a convergence theorem leading
to compactness, see~\cite{MS2}), the closure of $\X(A,J,p)$
contains points that lie on some cusp-curves representing the
class $A$ in a sense that we shall describe in a moment little
later. Therefore in order to compactify $\X(A,J,p)$ we must add
all simple cusp-curves in the class $A$ to the moduli space
$\M(A,J)$. The compactification is important because we want
$\X(A,J,p)$ to carry a fundamental homology class. We describe
below some features of a cusp-curve.

A cusp-curve $\s$ in $(M,\om)$, which represents the homology
class $A$, is a collection $\s=(\vt{\s}{N})$ of $J$-holomorphic
spheres $\s_i\hbox{:}\ \C P^1\lgra M$ such that $C_1\cup\cd\cup
C_N$ is a connected set, where $C_i=\s_i(\C P^1)$ and $A=A_1 + \cd
+ A_N$, $A_i$ being the homology class represented by $\s_i$. The
$\s_i$ are called the components of $\s$. A~cusp-curve $\s$ is
called simple if its components $\s_i$ are simple $J$-holomorphic
spheres such that $\s_i\ne \s_j\ox\ph$ for $i\ne j$ and any
$\ph\in G$. Any cusp-curve can be simplified to a simple
cusp-curve by replacing each multiply covered component by its
underlying simple curve. Of course this operation will change the
homology class $A$, but not the set of points that lie on the
curve. Also one can order the components of $\s$ so that
$C_1\cup\cd\cup C_k$ is connected for all $k\le N$. This means
that there exist integers $j_2,\ld, j_N$ with $1\le j_i\lt i$ such
that each $C_i$ must intersect some $C_{j_i}$, that is, there
exist $w_i$, $z_i\in \C P^1$ such that $\s_{j_i}(w_i)=\s_i(z_i)$.

A framing or intersection pattern $D$ of an ordered simple
cusp-curve
\begin{equation*}
\s=(\s_1,\ld, \s_N)
\end{equation*}
is a collection
\begin{equation*}
D=\{\vt{A}{N}, j_2, \ld, j_N\},
\end{equation*}
where $A_i=[C_i]\in H_i(M,\Z)$ and $j_i$ are integers with $1\le
j_i\le i-1$ chosen so that $C_i$ intersects $C_{j_i}$ (i.e.
$C_i\cap C_{j_i}\ne \emptyset$). Then $\om(A_i)\le \om(A)$, and so
there are only a finite number framings $D$ associated to $\s$.

For a fixed framing $D=\{\vt{A}{N}, j_2,\ld, j_N\}$, and a  $J\in \J(M,\om)$, let
\begin{equation*}
\M(\vt{A}{N}, J)= \M(A_1,J)\by \cd\by \M(A_N,J).
\end{equation*}
Let $\M(D,J)$  be the moduli space
\begin{equation*}
\M(D,J)\ci \M(\vt{A}{N}, J)\by (\C P^1)^{2N-2}
\end{equation*}
consisting of all $(\s, w,z)$ where $\s=(\vt{\s}{N})$, $\s_i\in
\M(A_i,J)$,
\begin{align*}
w=(w_2, \ld, w_N)\in (\C P^1)^{N-1}\quad \mb{ and }\quad
z=(z_2,\ld, z_N)\in (\C P^1)^{N-1},
\end{align*}
such that $\s$ is a simple cusp-curve with
$\s_{j_i}(w_i)=\s_i(z_i)$ for $i=2,\ld, N$.

For a generic $J$, $\M(D,J)$ will be an oriented manifold of
dimension
\begin{equation*}
2\sum_{j=1}^N\ c_1(A_j) + 2n+4(N-1).
\end{equation*}
The proof uses the extended evaluation map
\begin{equation*}
e_D\hbox{:}\ \M(\vt{A}{N}, J)\by (\C P^1)^{2N-2}\lgra M^{2N-2}
\end{equation*}
given by
\begin{equation*}
e_D(\s,w,z)= (\s_{j_2}(w_2),\s_2(z_2),\ld ,
\s_{j_N}(w_N),\s_N(z_N)).
\end{equation*}
The map $e_D$ is transversal to the multi-diagonal set
\begin{equation*}
\Delta_N=\{(x_2,y_2,\ld, x_N,y_N)\in M^{2N-2}|x_j=y_j\},
\end{equation*}
and therefore the inverse image $e_D^{-1}(\Delta_N)= \M(D,J)$ is a
manifold of the above dimension.

The group $G^N=G\by\cd\by G$ acts freely on $\M(D,J)$ by
\begin{align*}
\ph\cdot (\s_j,w_i,z_i)=(\s_j\ox \ph^{-1}_j, \ph_{j_i}(w_i),
\ph_i(z_i)),\ \ph=(\ph_1,\ld, \ph_N)\in G^N.
\end{align*}
The quotient space $\M(D,J)/G^N$ for a generic $J$ is a manifold
of dimension
\begin{equation*}
2c_1(A_1+\cd + A_N) +2n -2N- 4.
\end{equation*}
This is precisely our previous moduli space $\M(A,J)/G$ when $N=1$
and $A_1=A$.

Let $T$ denote a function $\{1,\ld, p\}\mapsto \{1,\ld, N\}$. This
function will indicate which of the $N$ components of
$C=C_1\cup\cd \cup C_N$ will be evaluated to get a point of
$M^p$.\break Define
\begin{equation*}
\W(D,T,J,p)=\M(D,J)\by_{G^N} (\C P^1)^p,
\end{equation*}
where the $j$th component of $\ph=(\vt{\ph}{N})\in G^N$ acts on
$\M(D,J)$ as above, and it acts on the $i$th factor of $(\C
P^1)^p$ if and only if $T(i)=j$. Then $\W(D,T,J,p)$ will be a
manifold of dimension
\begin{align*}
2\sum_{j=1}^N\ c_1(A_j) + 2n + 2p -2N -4.
\end{align*}
We have an evaluation map $e_{D,T}\hbox{:}\ \W(D,T,J,p)\lgra M^p$
defined by
\begin{equation*}
e_{D,T}(\s,w,z,\xi)=(\s_{T(1)}(\xi_1), \ld ,\s_{T(p)}(\xi_p)),
\end{equation*}
where $(\s,w,z)\in \M(D,J)$ and $\xi=(\vt{\xi}{p})\in (\C P^1)^p$.

We shall now choose $J$ suitably so that $\X(A,J,p)$ has a
fundamental homology class.

A manifold $(M,\om)$ is weakly monotone if every spherical
homology class $B\in H_2(M,\Z)$ with $\om(B)\gt 0$ and $c_1(B)\lt
0$ must satisfy the condition $c_1(B)\le 2-n$. Here $c_1$ is the
first Chern class of the complex bundle $(TM,J)$. This means that
there are no $J$-holomorphic spheres in homology classes with
negative first Chern number. The manifold $(M,\om)$ is monotone if
there is a $\l\gt 0$ such that $\om(B)=\l c_1(B)$ for every
spherical $B\in H_2(M,\Z)$. It can be shown that a monotone
manifold is weakly monotone, and conversely.

Let $R$ be a positive number. Then an $\om$-compatible almost
complex structure $J$ is called $R$-semi-positive if for every
$J$-holomorphic sphere $\s\hbox{:}\ \C P^1\lgra M$ with energy
$E(\s)\le R$ has Chern number $\int_{\C P^1} \s^*c_1\ge 0$. Let
$\J_+(M,\om,R)$ be the set of all
$\om$-compatible~$R$-semi-positive $J$. This set may be empty.
However, if $(M,\om)$ is a weakly monotone compact symplectic
manifold, then $\J_+(M,\om,R)$ is a path connected open dense set
for\break every $R$.

\setcounter{theore}{0}
\begin{theor}[\!]
Let $(M,\om)$ be a weakly monotone compact symplectic manifold{\rm
,} and $A\in H_2(M,\Z)$.
\begin{enumerate}
\leftskip .14pc
\renewcommand\labelenumi{{\rm (\alph{enumi})}}

\item For every $J\in \J(M,\om)$ there is a finite number of
evaluation maps
\begin{equation*}
e_{D,T}\hbox{\rm :}\ \W(D,T,J,p)\lgra M^p
\end{equation*}
such that
\begin{equation*}
\bigcap\ \ov{e_p(\W(A,J,p)-K)}\ci \bigcup_{D,T}\
e_{D,T}(\W(D,T,J,p)),
\end{equation*}
where the intersection is over all compact subsets $K$ in
$\W(A,J,p)${\rm ,} and the union is over all effective framings
$D$ and all functions $T\hbox{\rm :}\ \{1,\ld, p\}\lgra \{1,\ld,
N\}$.

\item There is a residual subset $\J_r$ in $\J(M,\om)$ such that
for every $J\in \J_r${\rm ,} the spaces $\W(D,T,J,p)$ are smooth
oriented $\s$-compact manifolds with
\begin{equation*}
\dim\ \W(D,T,J,p) = 2n + 2c_1(D) + 2p -2N-4.
\end{equation*}

\item Suppose that $A$ is not a multiple class $\l B$ where $\l\gt
1$ and $c_1(B)=0$. If $J\in \J_+(M,\om,R)\cap \J_r${\rm ,} then
\begin{equation*}
\dim\ \W(D,T,J,p)\le \dim\ \W(A,J,p) -2.
\end{equation*}
\end{enumerate}
\end{theor}\vspace{-.3pc}

Recall that a manifold is $\s$-compact if it is the union of a
countable family of compact sets.

The proof may be found in~\cite{MS2}. To understand the essence of
the theorem, we need to look at some facts about pseudo-cycles.

If a subset $X$ of a manifold $M$ is within the image of a smooth
map $g\hbox{:}\ V\lgra M$ defined on a smooth manifold $V$ of
dimension $k$, then $X$ is said to be of dimension at most $k$.
The boundary of the set $g(V)$, denoted by $g(V^\iy)$, is defined
to be the set
\begin{equation*}
g(V^\iy)= \bigcap\ \ov{g(V-K)},
\end{equation*}
where the intersection is over all compact subsets $K$ of $V$.
This is the set of all limit points of sequences $\{g(x_n)\}$
where $\{x_n\}$ has no convergent subsequence in $V$.

A $k$-pseudo-cycle in $M$ is a smooth map $f\hbox{:}\ W\lgra M$
defined on a smooth manifold $W$ of dimension $k$ such that $\dim\
f(W^\iy)\le k-1$. Two $k$-pseudo-cycles $f_0\hbox{:}\ W_0\lgra M$
and $f_1\hbox{:}\ W_1\lgra M$ are bordant if there is a
$(k+1)$-pseudo-cycle $F\hbox{:}\ W\lgra M$ with $\pl W=W_1-W_0$
such that
\begin{equation*}
F|W_0=f_0,\quad F|W_1=f_1,\quad \mb{ and }\quad \dim\ F(W^\iy)\le
k-1.
\end{equation*}

Every singular homology class $\a\in H_k(M)$ can be represented by
a $k$-pseudo-cycle $f\hbox{:}\ W\lgra M$. This can be seen in the
following way. First represent $\a$ by a map $f\hbox{:}\ X\lgra M$
defined on a finite oriented $k$-simplicial complex $X$ without
boundary so that $\a=f_*[X]$, where $[X]$ is the fundamental
class. Then approximate $f$ by a map which is smooth on each
simplex of $X$. Finally, consider the union of the $k$- and
$(k-1)$-faces of $X$ as a smooth manifold $W$ of dimension $k$ and
approximate $f$ by a map that is smooth across the
$(k-1)$-simplexes.

It also follows that bordant $k$-pseudo-cycles are in the same
homology class. However, two pseudo-cycles representing the same
homology class may not be bordant.

Theorem~5.1 says that the evaluation map $e_p\hbox{:}\
\W(A,J,p)\lgra M^p$ is a pseudo-cycle, and that the boundary of
its image can be covered by the sets
\begin{equation*}
e_{D,T}(\W(D,T,J,p)).
\end{equation*}
Therefore the image of $e_p$ carries a fundamental homology class.
It can be shown that this class is independent of the choice of
the point $z\in (\C P^1)^p$ and the almost complex structure $J$.

\section{Gromov--Witten invariants and quantum cohomology}

For the definition of the quantum cohomology, the symplectic
manifold $(M,\om)$ is required to satisfy the following (mutually
exclusive) conditions:
\begin{enumerate}
\leftskip .1pc
\renewcommand\labelenumi{{\rm (\alph{enumi})}}

\item $M$ is monotone, that is, $\langle\om,A\rangle=\l \langle c_1,A\rangle$ for every
$A\in \p_2(M)$, where $\l\gt 0$ and $c_1=c_1(TM,J)$.

\item $\langle c_1,A\rangle=0$ for every $A\in \p_2(M)$, or $\langle \om,A\rangle =0$ for
every $A\in \p_2(M)$.

\item The minimal Chern number $N$, defined by $\langle c_1,\p_2(M)\rangle =N\Z$
where $N\ge 0$, is greater than or equal to $n-2$.\vspace{-.52pc}
\end{enumerate}

It can be shown that a manifold $(M,\om)$ is weakly monotone if
and only if one of the above conditions is satisfied.

Let $(M,\om)$ be a weakly monotone compact symplectic manifold
with a fixed $A\in H_2(M,\Z)$. Then the $p$-fold evaluation map
\begin{equation*}
e_p\hbox{:}\ \W(A,J,p)\lgra M^p
\end{equation*}
represents a well-defined homology class in $M^p$, which is
independent of $J$.

If $A$ is a spherical homology class, and $p\ge 1$, then define a
homomorphism
\begin{equation*}
\Ph_{A,p}\hbox{:}\ H_d(M^p,\Z)\lgra\Z,
\end{equation*}
where $d=2np -\dim \W(A,J,p)$, in the following way. Let $\a\in
H_d(M^p,\Z)$ so that $\a=\a_1\by\cd\by\a_p$, where $\a_j\in
H_{d_j}(M,\Z)$ with $d_1 +\cd +d_p=d$. One can find a cycle
representing the homology class $\a$, which is denoted by the same
notation $\a$, such that it intersects the image $\X(A,J,p)$ of
the map $e_p$ transversely in a finite number of points. Then the
Gromov--Witten invariant $\Ph_A(\vt{\a}{p})$ is the intersection
number $e_p\cdot \a$, which is the number of intersection points
counted with signs according to their orientations. This is the
number of $J$-holomorphic spheres $\s$ in the homology class $A$
which intersect each of the cycles $\vt{\a}{p}$. If the dimension
condition for $d$ is not satisfied, then one sets
$\Ph_A(\vt{\a}{p})=0$.

The quantum cohomology is obtained by defining a quantum
deformation of the cup product on the cohomology of a symplectic
manifold $(M,\om)$. Before going into this, let us review the
ordinary cup product in singular cohomology.

We denote by $H^*(M)$ the free part of $H^*(M,\Z)$. We may
consider $H^*(M)$ as de Rham cohomology consisting of classes
which take integral values on all cycles:
\begin{equation*}
H^*(M)=H^*_{\rm DR}(M,\Z).
\end{equation*}
Next, we let $H_*(M)$ denote $H_*(M,\Z)/\mrm{Torsion}$. Then we
can identify $H^k(M)$ with $\mrm{Hom} (H_k(M),\Z)$ by the pairing
of $a\in H^k(M)$ and $\b\in H_k(M)$ given by
\begin{equation*}
a(\b)=\int_\b\ a.
\end{equation*}
In the same way, the intersection pairing $\a\cdot\b$ of $\a\in
H_{2n-k}(M)$ and $\b\in H_k(M)$ gives rise to the homomorphism
\begin{equation*}
\mrm{PD}\hbox{:}\ H_{2n-k}(M)\lgra H^k(M),
\end{equation*}
where $\mrm{PD}(\a)=a$ if
\begin{equation*}
a(\b)=\int_\b\ a= \a\cdot\b\quad \mb{ for } \ \b\in H_k(M).
\end{equation*}
The $\mrm{Poincar\acute{e}}$ duality theorem says that $\mrm{PD}$
is an isomorphism. Then the cup product $a\cup b\in H^{k+\ell}(M)$
of $a\in H^k(M)$ and $b\in H^\ell(M)$ is defined by the triple
intersection
\begin{equation*}
\int_\g\ a\cup b = \a\cdot\b\cdot\g,\quad \mb{ for } \ \g\in
H_{k+\ell}(M),
\end{equation*}
where $\a=\mrm{PD}^{-1}(a)\in H_{2n-k}(M)$ and
$\b=\mrm{PD}^{-1}(b)\in H_{2n-\ell}(M)$. This is well-defined,
because if the cycles representing $\a$ and $\b$ are in general
position, then they intersect a pseudo-cycle of codimension
$k+\ell$.

Next note that by our assumption $(M,\om)$ is monotone with
minimal Chern number $N\ge 2$.

The quantum multiplication $a*b$ of classes $a\in H^k(M)$ and
$b\in H^\ell(M)$ is defined as follows. Let $\a=\mrm{PD} (a)$ and
$\b=\mrm{PD} (b)$ denote the $\mrm{Poincar\acute{e}}$ duals of $a$
and $b$ so that $\deg(\a)=2n-k$ and $\deg(\b)=2n-\ell$. Then $a*b$
is the formal sum
\begin{equation*}
a*b = \sum_A\ (a*b)_A\cdot q^{c_1(A)/N},
\end{equation*}
where $q$ is an auxiliary variable supposed to be of degree $2N$,
and the cohomology class $(a*b)_A\in H^{k+\ell - 2c_1(A)}(M)$ is
defined in terms of the Gromov--Witten invariant $\Ph_A$ by
\begin{equation*}
\int_\g\ (a*b)_A =\Ph_A(\a, \b, \g),
\end{equation*}
for $\g\in H_{k+ \ell - 2c_1(A)}(M)$. Here $\a$, $\b$, $\g$
satisfy the following dimension condition required for the
definition of the invariant $\Ph_A$,
\begin{equation*}
2c_1(A) + \deg (\a) + \deg (\b) +\deg (\g) =4n.
\end{equation*}
The condition shows that $0\le c_1(A)\le 2n$, and therefore only
finitely many powers of $q$ occur in the above sum defining $a*b$.
Since $M$ is monotone, the classes $A$ which contribute to the
coefficient of $q^d$ satisfy $\om(A)=c_1(A)/N=d$, and therefore
only finitely many can be represented by $J$-holomorphic spheres.
Therefore the sum is finite. Since only nonnegative powers of $q$
occur in the sum, it follows that $a*b$ is an element of the group
\begin{equation*}
\widetilde{QH}^*(M)=H^*(M)\ot \Z[q],
\end{equation*}
where $\Z[q]$ is the polynomial ring in the variable $q$ of degree
$2N$. Then we get a multiplication by linear extension
\begin{equation*}
\widetilde{QH}^*(M)\ot \widetilde{QH}^*(M)\lgra
\widetilde{QH}^*(M).
\end{equation*}

The quantum cup product is skew-commutative in the sense that
\begin{equation*}
a*b=(-1)^{\deg a\cdot\deg b}b*a
\end{equation*}
for $a, b\in QH^*(M)$. Moreover, the product is distributive over
the sum, and associative. The skew-symmetry and the distributive
properties follow easily. But the associative property is a bit
complicated and depends on a certain gluing argument for
$J$-holomorphic spheres.

The quantum cohomology $\widetilde{QH}^k(M)$ vanish for $k\le 0$,
and are periodic with period $2N$ for $k\ge 2n$. To get the full
periodicity, we consider the group
\begin{equation*}
QH^*(M) =H^*(M)\ot \Z[q, q^{-1}],
\end{equation*}
where $\Z[q, q^{-1}]$ is the ring of Laurent polynomials, which
consists of polynomials in the variables $q$, $q^{-1}$ with the
obvious relation $q\cdot q^{-1}=1$. With this definition $QH^k(M)$
is non-zero for positive and negative values of $k$, and there is
a natural isomorphism
\begin{equation*}
QH^k(M)\lgra QH^{k+2N}(M)
\end{equation*}
given by multiplication with $q$, for every $k\in \Z$.

If $A=0$, then all $J$-holomorphic spheres in the class $A$ are
constant. It follows then that $\Ph_A(\a,\b,\g)$ is just the usual
triple intersection $\a\cdot\b\cdot\g$. Since $\om(A)> 0$ for all
other $A$ which have the $J$-holomorphic representatives, the
constant term in $a*b$ is just the usual cup product.

The product in $QH^*(M)$ is also distributive over the sum, and
skew-commutative. It commutes with the action of $\Z[q,q^{-1}]$.
If $a\in H^0(M)$ or $H^1(M)$, then $a*b=a\cup b$ for all $b\in
H^*(M)$. Also the canonical generator ${\bf 1}\in H^0(M)$ is the
unit in quantum cohomology.

As an example, let $M$ be the complex projective $n$-space $\C
P^n$ with the standard $\mrm{K\ddot{a}hler}$ form. Let $L$ be the
standard generator of $H_2(\C P^1)$ represented by the line $\C
P^1$. Then the first Chern class of $\C P^n$ is given by
$c_1(L)=n+1$. Therefore, by the dimension condition, the invariant
$\Ph_{mL}(\a,\b,\g)$ is non-zero only when $m=0$ and $1$. Clearly,
the case $m=0$ corresponds to constant curves, and gives the
ordinary cup product. Since the minimal Chern number is $N=n+1$,
the quantum cohomology groups are given by $QH^k(M)\simeq \Z$ when
$k$ is even, and $QH^k(M)=0$ when $k$ is odd.

Next, let $a\in H^\ell(M)$ and $b\in H^k(M)$. If $\ell+k\le 2n$,
then the quantum cup product is the same as the ordinary cup
product $a*b=a\cup b$. So consider the case when $a$ is the
standard generator $p$ of $H^2(M)$ defined by $p(L)=1$, and
$b=p^n\in H^{2n}(M)$. Then the quantum cup product $p*p^n$ is the
generator $q$ of $QH^{2n+2}(M)$, because
\begin{equation*}
\int_{pt} (p*p^n)_L=\Ph_L([\C P^{n-1}], pt, pt)=1,
\end{equation*}
where $[\C P^{n-1}]=PD(p)$ and $pt=PD(p^n)$, and all other classes
$(p*p^n)_A$ vanish. Therefore the quantum cohomology of $\C P^n$
is given by
\begin{equation*}
\widetilde{QH}^*(\C P^n)=\f{\Z[p,q]}{(p^{n+1}=q)}.
\end{equation*}

\section*{Acknowledgement}
The author is grateful to the referee for pointing out certain
inaccuracies.

\end{document}